\documentclass[12pt,plain]{article}
\usepackage{listings}
\usepackage[makeroom]{cancel}
\usepackage{fancyhdr}
\usepackage{verbatim}
\usepackage{indentfirst}
\usepackage{graphicx}
\usepackage{color}
\usepackage{newlfont}
\usepackage{amssymb}
\usepackage[intlimits]{amsmath}
\usepackage{latexsym}
\usepackage{amsthm}
\usepackage{amscd}
\usepackage{fullpage}
\usepackage{setspace}
\usepackage{multirow}
\usepackage{hyperref}
\usepackage[dvips, hmargin=2.1cm,vmargin=3cm]{geometry}
\usepackage{multirow}
\usepackage{array}

\newcommand{\N}{\mathbb{N}}

\newcommand{\R}{\mathbb{R}}

\newcommand{\be}{\begin{equation}}
\newcommand{\bey}{\begin{eqnarray}}
\newcommand{\ee}{\end{equation}}
\newcommand{\eey}{\end{eqnarray}}
\newcommand{\ba}{\begin{array}}
\newcommand{\ea}{\end{array}}

\theoremstyle{plain}
\newtheorem{theorem}{Theorem}[section]

\newtheorem{cor}[theorem]{Corollary}
\theoremstyle{definition}

\newtheorem{remark}[theorem]{Remark}

\newtheorem{conj}[theorem]{Conjecture}

\input xy
\xyoption{all}

\begin{document}
\title{On two conjectures for M$\&$m sequences}
\author{F. Cellarosi$^*$, S. Munday$^\dag$}

\maketitle

\let\oldthefootnote\thefootnote
\renewcommand{\thefootnote}{\fnsymbol{footnote}}
\footnotetext[1]{University of Illinois at Urbana-Champaign. Department of Mathematics, Altgeld Hall, 1409 W. Green Street, Urbana, IL 61801,  United States. \texttt{fcellaro@illinois.edu}}
\footnotetext[2]{University of York. Department of Mathematics, York YO10 5DD, United Kingdom. \texttt{sam584@york.ac.uk}}
\let\thefootnote\oldthefootnote

\begin{abstract}
In this paper, the recently introduced M$\&$m sequences and associated mean-median map are studied. These sequences are built by adding new points to a set of real numbers by balancing the mean of the new set with the median of the original. This process, although seemingly simple, gives rise to complicated dynamics. The main result is that two conjectures put forward by Chamberland and Martelli are shown to be true for a subset of possible starting conditions.
\end{abstract}

\begin{keywords}
M\&m sequences, Mean-median map, strong terminating conjecture, continuity conjecture
 {\bf MSC}: 97F40, 11B75, 37E15
\end{keywords}

\section{Introduction}
The aim in this paper is to continue the exploration of mean-median sequences and the associated mean-median map. The class of mean-median sequences, the generation of which we shortly describe, was introduced by Schultz and Shiflett \cite{SS} and further analysed by Chamberland and Martelli \cite{CM}, who introduced the mean-median map to aid their investigations, and also by Bonchev Bonchev \cite{BB}. The interest in the mean-median map lies in the fact that it is a very good example of a simple process yielding extremely complicated dynamics. Let us now give the details of the construction.

The \emph{median} of a finite set of real numbers $\{\xi_1,\xi_2,\ldots,\xi_N\}$, with
$\xi_{i_1}\leq\xi_{i_2}\leq\ldots\leq \xi_{i_N}$, is given by
\be\mathrm{med}(\xi_1,\ldots,\xi_N)=\begin{cases}\xi_{i_{(N+1)/2}},&\mbox{$N$ odd;}\\&\\ \displaystyle\frac{\xi_{i_{N/2}}+\xi_{i_{N/2+1}}}{2},&\mbox{$N$ even.}\end{cases}\nonumber\ee
Given three real numbers $a,b,c$, we want to add a fourth number $x_4$  to the list so that the mean of the four numbers equals the median of the three. That is, $x_4$ is defined to be the unique solution to the equation
\be\frac{a+b+c+x_4}{4}=\mathrm{med}(a,b,c).\nonumber\ee
Iteratively, we let $x_n$ be the unique solution of the equation
\be\frac{a+b+c+x_4+\ldots+x_n}{n}=\mathrm{med}(a,b,c,x_4,\ldots,x_{n-1}).\label{def-process}\ee

Schultz and Shiflett gave the name {\em M$\&$m sequences}, for mean and median, to the sequences $(x_4, x_5, \ldots)$ generated as in \eqref{def-process}. Their paper begins with the observation that no matter which set of three numbers $\{a, b, c\}$ they picked to begin the process, the resulting M$\&$m sequence was eventually constant, that is, they found a  $k\in \N$ such that $x_n=x_k$ for all $n\geq k$. They called an M$\&$m sequence with this property {\em stable} and conjectured that all such sequences were stable. (It is perhaps worth observing here that while it is certainly possible to find examples which stabilise very quickly, it is also possible to find examples which take an extremely long time to stabilise, as will be apparent in our main results below.) They also pointed out that it is enough to consider sequences generated from sets $\{a, b, c\}$ with $a<b<c$, since if all three are equal the sequence becomes stable immediately, and if two of the starting values are equal, the sequence becomes stable after five iterations.

By applying an affine transformation, we can always reduce the case of arbitrary $a<b<c$ to $0<x<1$ (this is  shown in \cite{CM} and \cite{SS}, although in \cite{SS} they choose the different normalisation $0<x<x+1$). So, from this point on, let us suppose that we start with $0<x<1$ and consider the M$\&$m sequence $(x_4,x_5,\ldots)$ generated as in \eqref{def-process}. We can also formulate questions about the sequence $(x_4,x_5,\ldots)$ in terms of the sequence $(m_4,m_5,\ldots)$, where $m_n=\mathrm{med}(0,x,1,x_4,\ldots,x_n)$. One can easily check that the sequence of medians $(m_n)_{n\geq 4}$ is monotone (see Theorem 2.1 in \cite{CM}).  Schultz and Shiflett's stability conjecture was reformulated in the following way in \cite{CM}.

\begin{conj}[Strong Terminating Conjecture, \cite{CM}]\label{strong-terminating-conj}
For every $0<x<1$, there exists an integer $k$ such that $x_n=x_k$ for all $n>k$.
 The minimum such $k$ is denoted by $L(x)$ and is referred to as the \emph{length} of the M\&m sequence. We can reformulate the conjecture by saying that for every initial choice of $0<x<1$, we have that  $L(x)<\infty$ or, equivalently, 
 that the sequence of medians eventually becomes constant.
\end{conj}

We can now define the {\em mean-median map}, as in \cite{CM}. Let  $m:\R\to\R$ be defined to be the function that assigns to each $x$ with $0<x<1$ the limit of the monotone sequence of medians $(m_n)_{n\geq 4}$ associated to it, provided that this limit exists. In other words, if the strong terminating conjecture is true, $m(x)$ is equal to the constant value achieved by the sequence of medians associated to $x$ in a finite number of steps. The following conjecture about the function $m$ has also been made.

\begin{conj}[Continuity Conjecture, \cite{CM}]\label{continuity-conj}
The function $x\mapsto m(x)$ is continuous.
\end{conj}

If the Strong Terminating Conjecture  and the Continuity Conjecture turn out to be true, then it follows that the function $m$ is piecewise affine, with rational corners.

As shown in Theorems 3.1 and 3.2 in \cite{CM}, the problem has two symmetries that allow us to restrict our study to the subinterval $\frac{1}{2}\leq x\leq \frac{2}{3}$. These are:
\begin{align}
m(1-x)&=1-m(x)&0\leq x\leq 1\label{symmetry1}\\
m(x)&=(3x-1)m\!\left(\frac{x}{3x-1}\right)&\frac{1}{2}\leq x\leq 1\nonumber
\end{align}

Our main result is the following
\begin{theorem}\label{main-thm}
For every rational $\frac{p}{q}\in\left[\frac{1}{2},\frac{2}{3}\right]$ with $2\leq q\leq 18$, conjectures \ref{strong-terminating-conj} and \ref{continuity-conj} are true in a neighborhood of $\frac{p}{q}$. Moreover, near $\frac{1}{2}$,
\begin{align}
m(x)&=\begin{cases}
\frac{333}{8}x-\frac{325}{16},&\frac{1}{2}\leq x\leq \frac{2911001}{5684610}\\
-\frac{2840973}{32}x+\frac{2909701}{64},&\frac{2911001}{5684610}\leq x\leq \frac{339}{662}\\
\frac{2842297}{32}x-\frac{2910929}{64},&\frac{339}{662}\leq x\leq \frac{56346353}{110033282}\\
\frac{333}{8}x-\frac{325}{16},&\frac{56346353}{110033282}\leq x\leq \frac{90028408624696264845}{175783207694477729162}\\
-\frac{87891603847238182597}{16384}x+\frac{90028408624695599245}{32768},&\frac{90028408624696264845}{175783207694477729162}\leq x\leq \frac{120193266020529}{234681009773714}\\
\frac{87885645878322471973}{16384}x-\frac{90022305806717048045}{32768},&\frac{120193266020529}{234681009773714}\leq x \leq \frac{90022305806716382445}{175771291756643579978}.
\end{cases}
\label{near1/2}
\end{align}
Near $\frac{2}{3}$,
\begin{align}
m(x)&=\begin{cases}
-\frac{75675009}{256}x+\frac{25065033}{128},&\frac{50110610}{75646209}\leq x\leq \frac{626}{945}\\
\frac{75647841}{256}x-\frac{25055657}{128},&\frac{626}{945}\leq x\leq \frac{50130770}{75676641}\\
-\frac{225}{2}x+76,&\frac{50130770}{75676641}\leq x\leq \frac{2}{3}.
\end{cases}
\label{near2/3}
\end{align}
Similar explicit formul\ae\, are available near all the other rational points, see Appendix A.
\end{theorem}
\begin{cor}
The function $x\mapsto m(x)$ has corners at the points
\begin{align}
\frac{9}{17}, \frac{8}{15}, \frac{7}{13}, \frac{6}{11}, \frac{5}{9}, \frac{9}{16}, \frac{4}{7}, \frac{7}{12}, \frac{10}{17}, \frac{3}{5}, \frac{11}{18}, \frac{8}{13}, \frac{5}{8}, \frac{7}{11}, \frac{9}{14}, \frac{11}{17}, \frac{2}{3}.\nonumber
\end{align}
\end{cor}
%

\begin{remark}
Our method is also able to establish whether a rational point is \emph{not} a corner for the function $m(x)$, like $x=\frac{1}{2}$ because of \eqref{symmetry1}.
For example, $\frac{10}{19}\in(\frac{841}{1598},\frac{639}{1214})$ and on this interval conjectures \ref{strong-terminating-conj} and \ref{continuity-conj} are valid with $L(x)=47$ and $m(x)=\frac{141}{4}x-\frac{137}{8}$.
\end{remark}

Theorem \ref{main-thm} improves the results of Chamberland and Martelli, who proved in \cite{CM} that $m(x)=\frac{333}{8}x-\frac{325}{16}$
for 
$\frac{1}{2}\leq x\leq \frac{337}{666}$
and claimed
that $m(x)=\frac{225}{2}x+76$ for points  $x$ sufficiently  close to $\frac{2}{3}$. The proof of Theorem \ref{main-thm} is computer-assisted, in that it uses an algorithm to find explicit neighbourhoods of rational points where the two conjectures hold true, and derive the exact formula for $m(x)$ in these neighbourhoods. In  Section \ref{section-method} we explain our algorithm. Section \ref{section-full-code} includes our implementation of the method, along with the derivation of \eqref{near1/2} and some discussion on the combinatorial features of the problem.

\section*{Acknowledgments}
The authors would like to thank Franco Vivaldi for bringing the subject of mean-median sequences to their attention. They are also grateful to the Max Planck Institute for Mathematics in Bonn for providing access to their computers during the course of this investigations and a very congenial working environment for the first author during his visit there in June 2014. The first author gratefully acknowledges the financial support of the AMS-Simons travel grant (2013-2014) and the NSF grant DMS-1363227.

\section{Our method}\label{section-method}
Let us describe our proof of Theorem \ref{main-thm}.
Given a point $\frac{1}{2}\leq  x_0\leq \frac{2}{3}$ we find adjacent closed intervals
\be 
I^{(0)},I^{(1)},I^{(2)},\ldots\nonumber\ee
with the following properties:
\begin{itemize}
\item[(i)] $x_0\in I^{(0)}$.
\item[(ii)] For every $i$, the interval $I^{(i)}$ is a finite disjoint union of subintervals $I^{(i)}_k$, $1\leq k\leq K_i$ for some positive $K_i$. The intervals $I^{(i)}_k$ can be open at both ends, closed at both ends, or closed at only one end, and can even consist of a single point.
\item[(iii)] For every $x\in I^{(i)}_k$ the length of the M\&m sequence associated to $x$ is constant. This length is different from the length(s) in adjacent interval(s) $I^{(i)}_{k}$.
\item[(iv)] The function $x\mapsto m(x)$ is affine on each $I^{(i)}$, with corners at the endpoints of $I^{(i)}$.
\end{itemize}
Our algorithm begins by constructing the interval $I^{(0)}_s$ containing $x_0$ in its closure. Then the adjacent intervals $I^{(0)}_{s+1}, I^{(0)}_{s+2}, \ldots, I^{(0)}_{K_0}$ are found sequentially, left to right. The intervals on the opposite side $I^{(0)}_1,\ldots,I^{(0)}_{s-1}$ can be computed independently, right to left. This feature can be exploited by a parallel implementation of our algorithm. 
On each of the intervals $I^{(0)}_k$ an explicit formula for $L(x)$ and $m(x)$ is found by solving a finite number or linear inequalities with rational coefficients. When transitioning from $I^{(0)}_k$ to $I^{(0)}_{k+1}$ the value of $L$ increases or decreases, but the affine formula for $m(x)$ remains the same. If the formula for $m(x)$ changes, then $K_0=k$ and the new interval is labeled $I^{(1)}_1$. Then we continue onwards with $I^{(1)}_2, I^{(1)}_{3},\ldots,I^{(1)}_{K_1}$, etc.
Similarly, we can determine analogous adjacent intervals $\ldots, I^{(-3)},I^{(-2)}, I^{(-1)}$ to the left of $I^{(0)}$.

Let us now describe the algorithm that produces the intervals $I^{(i)}_k$ and computes the functions $L(x)$ and $m(x)$ restricted to these intervals.
 The code we include is written using  \emph{Mathematica} language. Let us stress that the results obtained with our method are \emph{exact} and involve no numerical approximation. Only exact arithmetic and symbolic manipulations are used. The only limitation of our method is given by time and memory constraints related to the implementation of the algorithm.
Here is the function we use to compute the median of a vector.
{\small
\begin{verbatim}
median[u_] := If[OddQ[Length[u]],
                 u[[Ceiling[Length[u]/2]]],
                 (u[[Length[u]/2]] + u[[Length[u]/2 + 1]])/2];
\end{verbatim}
}
Notice that, in contrast with the built-in \texttt{Median[]}, the function above also works for chains of inequalities, which will be needed later. For example, \texttt{median[2 < a < b < Pi < c < 6]} yields \texttt{(b + Pi)/2}.

The simple routine that we use to find $L(x)$ and $m(x)$ for a given number $x$ is as follows. Without loss of generality we can assume that $x$ is rational so that exact arithmetic can be used. We start  with the list $(0,x,1)$. Let $M=x$ (the median of the list) and $S=x+1$ (the sum of the numbers in the list). We construct $x_4= 4 M-S$ as in \eqref{def-process}, and update the median $M$ and the sum $S$. Then proceed with the computation of $x_5$, etc.
{\small
\begin{verbatim}
threshold=10000;
Lm[x_] := (listx = {0, x, 1}; M = x; listM = {x}; S = x + 1;
   numstepsSTART = 0;
   Do[newx = FullSimplify[j M - S];
    S = S + newx;
    AppendTo[listx, newx];
    If[M == newx, Return[{j, Last[listM]}]; Break[]];
    M = Median[listx];
    AppendTo[listM, M], {j, 4, threshold}]);
\end{verbatim}
}

For example, if $x=\frac{7}{12}$, we have \texttt{listx = \{0, 7/12, 1, 3/4, 1, 7/6, 13/8, 15/8, 1\}}, \texttt{listM = \{7/12, 2/3, 3/4, 7/8, 1, 1\}} and the algorithm produces \texttt{Lm[7/12] = \{9,1\}}.

To compute the endpoints of the intervals $I^{(i)}_k$,  we start by finding the interval $I^{(0)}_s$ which includes $x_0$.
Besides $x_0\in\mathbb{Q}$, we 
also consider a small parameter $\epsilon>0$. This parameter can also be taken to be rational to allow exact arithmetic on a computer.
As we shall see, the interval $I^{(0)}_s$ will not depend on $\epsilon$, provided $\epsilon$ is small enough. The algorithm will also reduce $\epsilon$ if necessary.
%

First, we run the routine described above for the number $x'=x_0+\epsilon_0$  and obtain the corresponding M\&m sequence $\texttt{listx}=\{0,x',1,x_4,x_5,\ldots, x_{L(x')}\}$ with $m(x')=x_{L(x')}$.   It is known that, in order for  $m(x)=\bar x$ for some $x$, it is necessary (but not sufficient) that there exists $4\leq l<L(x)$ such that $x_l=\bar x$; see \cite{BB}. This means that \texttt{listx} includes $m(x')$ (and possibly other elements) more than once.

Then we use \texttt{listx} to create a \emph{driving list}, consisting of the indices of the elements of the M\&m sequence of $x'$, after sorting it in increasing order. In other words, the driving list is a permutation of $\{1,\ldots,L(x')\}$ describing the ordering of the trajectory of $x'$.
\begin{verbatim}
sortedlist = Sort[listx, Less];
drivinglist = DeleteDuplicates[Flatten[Table[Position[listx, sortedlist[[k]]],
  {k,1,Length[sortedlist]}]]];
\end{verbatim}

The next step of the algorithm finds all the $x$'s that yield the same driving list as $x'$. This set is obtained by solving a system of at most $L(x')$ linear inequalities with rational coefficients and consists of an interval containing $x'$ (but not necessarily $x_0$).
We start with the inequality $0<x<1$, corresponding to the permutation $(1,2,3)$ in the driving list. The median of this inequality is $M=x$, and the sum of the terms is $S=x+1$. The next point in the M\&m sequence, $x_4$, is computed as $x_4=4 M-S=-1+3x$. We now insert $x_4$ in the chain of inequalities, at the position prescribed by \texttt{drivinglist}. For example, if the numbers $\{1,2,3,4\}$ appear in \texttt{drivinglist} ordered as $(1,2,4,3)$. then the new inequality is $0<x<-1+3 x<1$. Then the new median $M$ and the sum $S$ are computed for the list of inequalities. In the example, we have $M=\frac{1}{2}(-1+4 x)$ and
 $S=4x$. The next point in the trajectory is now $x_5=5 M-S$ (equal to $-\frac{5}{2}+6 x$ in the example) and this inserted in the previous chain of inequalities according to the permutation of $\{1,2,3,4,5\}$ contained in \texttt{drivinglist}.
In the case of $(1,2,4,5,3)$ we get $0<x<-1+3 x<-\frac{5}{2}+6 x<1$. We continue this procedure until $x_{L(x')}$ has been placed in the chain of inequalities according to the driving list.  Here is the code to obtain the chain of inequalities from the driving list.
{\small
\begin{verbatim}
Clear[x];
permutation = {1, 2, 3};
inequalities = 0 < x < 1;
S = x + 1;
M = x;
listM = { };
Do[
   M = median[inequalities];
   AppendTo[listM, M]'
   S = Simplify[S + newx];
   newx = Simplify[n M - S];

   possiblepermutations = Table[Insert[permutation, n, t], {t, 1, n}];
   pos = Position[possiblepermutations, Select[drivinglist, # <= n &]][[1]];

   permutation = Insert[permutation, n, pos];
   inequalities = Insert[inequalities, newx, pos],
{n, 4, L}];
finalx=Last[listM];
\end{verbatim}
}
Recall that some elements in the trajectory $\texttt{listx}$  repeat 
and, since we are considering strict inequalities,  we 
remove the duplicates from the chain. 
Then, we simplify the chain of inequalities to get an open interval.
{\small
\begin{verbatim}
todelete = Position[Table[inequalities[[i]] == inequalities[[i + 1]],
                    {i, 1, L - 1}], True];
inequalities = Reduce[Delete[inequalities, todelete]];
 \end{verbatim}
}
If \texttt{inequalities} does not contain $x_0$ in its closure, then we repeat the procedure so far with a smaller $\epsilon$. Otherwise, we proceed as follows. All the points $x$ in the interior of the interval \texttt{inequalities} are have, by construction, $L(x)=$ \texttt{L} and $m(x)=$ \texttt{finalx}, which is an affine function of $x$ with rational coefficients.
The endpoints of the interval given by \texttt{inequalities} are checked separately,  to see if the values of $L(x)$ and $m(x)$ extend to the endpoints. If so, the strict inequalities are replaced by loose inequalities.
The new interval (along with the value of $L(x)$ and $m(x)$ on that interval) is then output to the list \texttt{results}. If $L(x)$ is different compared to the previous interval $I^{(i)}_k$, then the new interval is $I^{(i)}_{k+1}$. If the affine formula for $m(x)$ is different, then the new interval is $I^{(i+1)}_1$.

\section{The full code and an example}\label{section-full-code}
Here we include the full \emph{Mathematica} code for the algorithm we just described. We use the input $x_0=\frac{1}{2}$ and $\epsilon=\frac{1}{100000}$ as an example.

{\small
\begin{verbatim}
results = {};
endpointssubintervals = {};
endpointsintervals = {};
x0 = 1/2;
eps = 1/100000;
threshold = 2000;
howmanycycles = 500;
median[u_] := If[OddQ[Length[u]], u[[Ceiling[Length[u]/2]]],
   (u[[Length[u]/2]] + u[[Length[u]/2 + 1]])/2];
Lm[x_] := (listx = {0, x, 1}; M = x; listM = {x}; S = x + 1;
   Do[newx = FullSimplify[j M - S];
    S = S + newx;
    AppendTo[listx, newx];
    If[M == newx, Return[{j, Last[listM]}]; Break[]];
    M = Median[listx];
    AppendTo[listM, M], {j, 4, threshold}]);
Do[
 Clear[x];
 x = x0 + eps;
 L = Lm[x][[1]];
 sortedlist = Sort[listx, Less];
 drivinglist = DeleteDuplicates[Flatten[Table[Position[listx, sortedlist[[k]]],
 {k, 1, Length[sortedlist]}]]];
 Clear[x];
 permutation = {1, 2, 3};
 inequalities = 0 < x < 1;
 S = x + 1;
 newx = 0;
 M = x;
 listM = {};
 Do[
  M = median[inequalities];
  AppendTo[listM, M];
  S = Simplify[S + newx];
  newx = Simplify[n M - S];
  possiblepermutations = Table[Insert[permutation, n, t], {t, 1, n}];
  pos = Position[possiblepermutations, Select[drivinglist, # <= n &]][[1]];
  permutation = Insert[permutation, n, pos];
  inequalities = Insert[inequalities, newx, pos];,
  {n, 4, L}];
 finalx = Last[listM];
 todelete =
  Position[Table[
    inequalities[[i]] == inequalities[[i + 1]], {i, 1,
     Length[inequalities] - 1}], True];
 inequalities = Reduce[Delete[inequalities, todelete]];
 If[Length[inequalities] == 5 &&
   inequalities[[1]] <= x0 <= inequalities[[5]],
  If[{L, finalx /. x -> inequalities[[1]]} == Lm[inequalities[[1]]],
   inequalities[[2]] = LessEqual];
  If[{L, finalx /. x -> inequalities[[5]]} == Lm[inequalities[[5]]],
   inequalities[[4]] = LessEqual];
  toadd = {inequalities, L, finalx};
  If[Length[results] == 0,
   Print["The first interval starts at ", inequalities[[1]]];
   Print["In this interval L = " <> ToString[L] <> " and m(x) = ",
    finalx]; AppendTo[endpointssubintervals, inequalities[[1]]]];
  If[Length[results] > 0,
   If[Last[results][[2]] != L,
    Print["An interval with L = " <> ToString[L] <>
      " was found, starting at ", inequalities[[1]]];
    AppendTo[endpointssubintervals, inequalities[[1]]]];
   If[Last[results][[3]] =!= finalx,
    Print["An interval with different m(x) was found, starting at ",
     inequalities[[1]]];
    Print["In this interval L = " <> ToString[L] <> " and m(x) = ",
     finalx]; AppendTo[endpointsintervals, inequalities[[1]]];
    AppendTo[endpointssubintervals, inequalities[[1]]]]];
  AppendTo[results, toadd];
  x0 = inequalities[[5]],
  eps = eps/10],
 {howmanycycles}]
 \end{verbatim}
 }
The output for 500 cycles required 134.6 seconds on one of the authors' computer. We have
{\small
\begin{verbatim}
The first interval starts at 1/2
In this interval L = 73 and m(x) = 1/16 (-325+666 x)
An interval with L = 75 was found, starting at 341/666
An interval with L = 77 was found, starting at 24073/47010
An interval with L = 79 was found, starting at 24751/48334
An interval with L = 81 was found, starting at 24749/48330
An interval with L = 83 was found, starting at 784/1531
An interval with L = 85 was found, starting at 76279/148958
An interval with L = 87 was found, starting at 76957/150282
An interval with L = 89 was found, starting at 3263/6372
An interval with L = 91 was found, starting at 52547/102614
An interval with L = 93 was found, starting at 133909/261498
An interval with L = 95 was found, starting at 134587/262822
An interval with L = 97 was found, starting at 82379/160870
An interval with L = 99 was found, starting at 41359/80766
An interval with L = 101 was found, starting at 196963/384630
An interval with L = 103 was found, starting at 197641/385954
An interval with L = 105 was found, starting at 57631/112542
An interval with L = 107 was found, starting at 115601/225746
\end{verbatim}}
This means that $I^{(0)}_1=\{\frac{1}{2}\}$, $I^{(0)}_2=(\frac{1}{2},\frac{341}{666}]$, $I^{(0)}_3=(\frac{341}{666},\frac{24073}{47010}]$, $I^{(0)}_4=(\frac{24073}{47010},\frac{24751}{48334}]$, $I^{(0)}_{5}=(\frac{24751}{48334},\frac{24749}{48330}]$, and so on. The information about the endpoints is contained in the output \texttt{results} as we described.

To discover the first corner where the function $m(x)$ changes from $\frac{333}{8}x-\frac{325}{16}$ to $-\frac{2840973}{32}x+\frac{2909701}{64}$ is the same as finding the right endpoint for $I^{(0)}$. In order to achieve this, one has either to run the algorithm for much longer, or start with $x_0$ already close to the corner. For $x_0=\frac{5756575}{11241454}$ (the right endpoint of $I^{(0)}_{98}$ to the right of $\frac{1}{2}$) we get
{\small
\begin{verbatim}
The first interval starts at 5756575/11241454
In this interval L = 269 and m(x) = 1/16 (-325+666 x)
An interval with L = 271 was found, starting at 5757253/11242778
An interval with different m(x) was found, starting at 2911001/5684610
In this interval L = 271 and m(x) = 1/64 (2909701-5681946 x)
An interval with different m(x) was found, starting at 339/662
In this interval L = 271 and m(x) = 1/64 (-2910929+5684594 x)
An interval with different m(x) was found, starting at 2909629/5681930
In this interval L = 271 and m(x) = 1/16 (-325+666 x)
\end{verbatim}}
This means that $I^{(0)}=[\frac{1}{2},\frac{2911001}{5684610}]$, $I^{(1)}=[\frac{2911001}{5684610},\frac{339}{662}]$, $I^{(2)}=[\frac{339}{662},\frac{2909629}{5681930}]$, etc. See \eqref{near1/2}. All the other formul\ae\, in Theorem \ref{main-thm} are obtained in the same way. In the case of \eqref{near1/2} and \eqref{near2/3} the endpoints provided are corners of $m(x)$.

\begin{figure}[htbp]
\begin{center}
\includegraphics[angle=-90,width=5in]{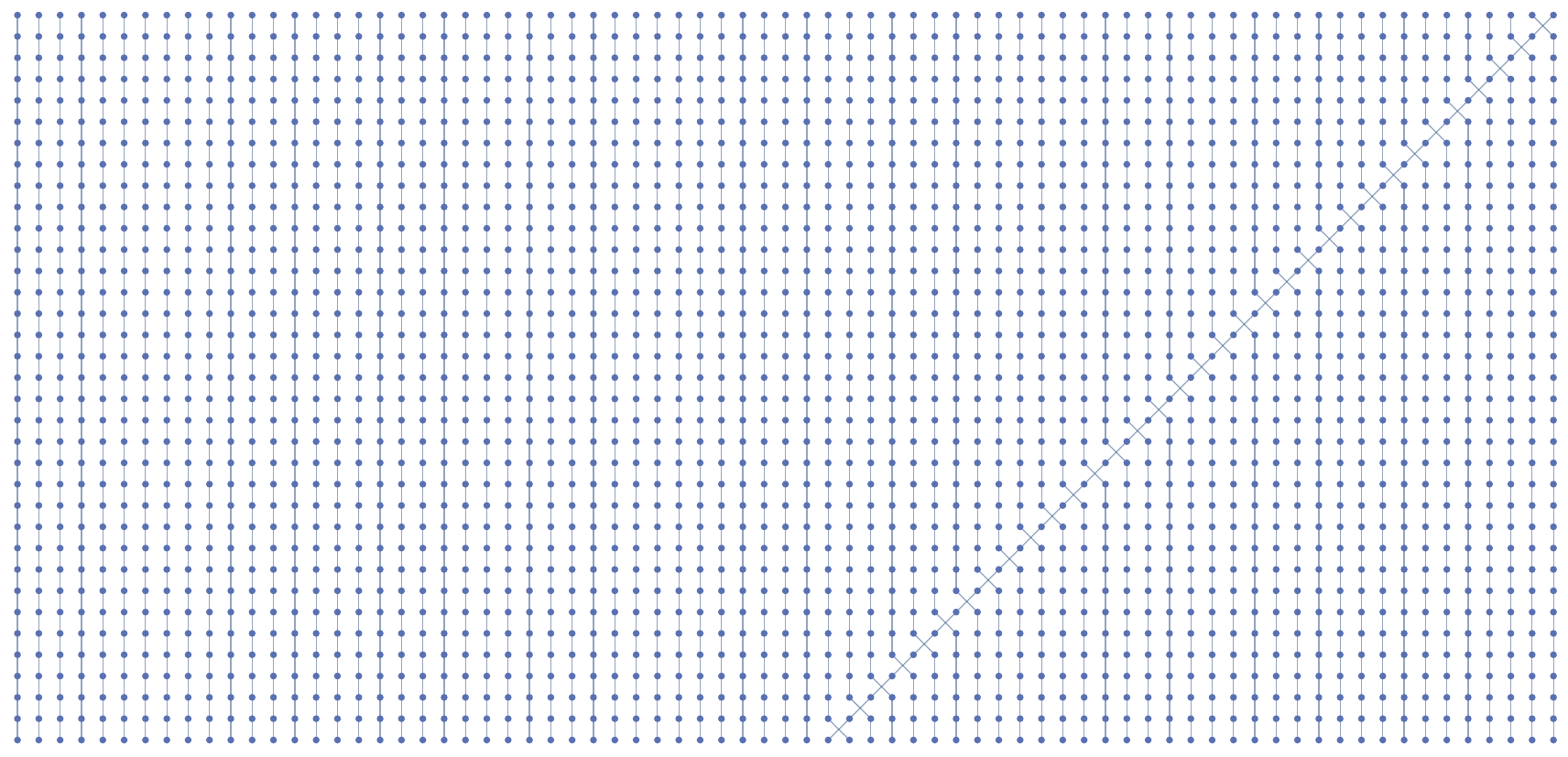}
\includegraphics[angle=-90,width=5.10in]{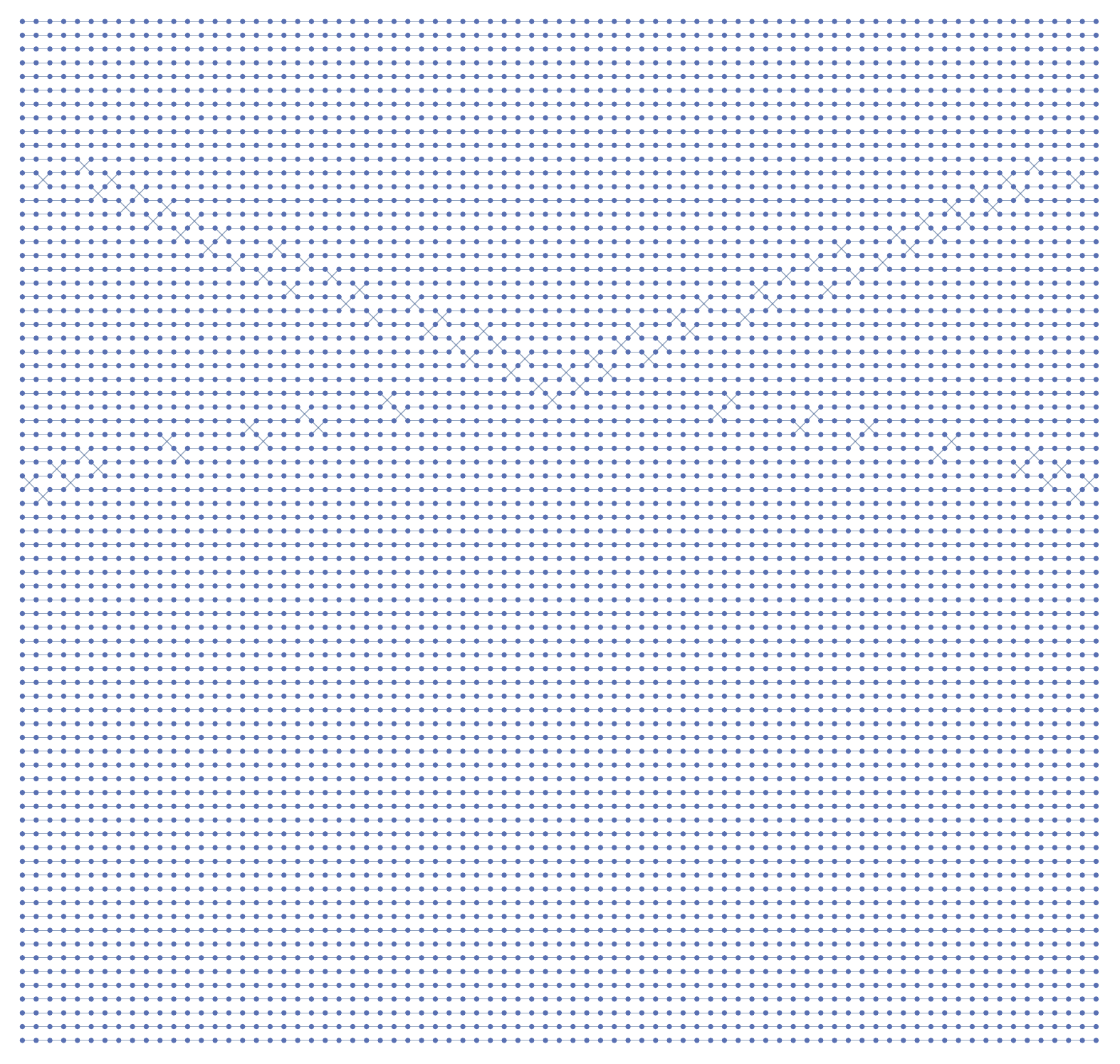}
\caption{Changes in the combinatorics of the M\&m sequence for $x\in I^{(0)}_2=(\frac{1}{2},\frac{341}{666}]$ (top) and for $x\in I^{(0)}_3=(\frac{341}{666},\frac{24073}{47010}]$ (bottom).}
\label{figure-permutations1}
\end{center}
\end{figure}

Another way to visualise the algorithm we described in Section \ref{section-method} is to see how the driving lists change within each interval $I^{(i)}_k$, that is,  how the combinatorics of the M\&m sequence changes when $x$ runs through the interior of the interval $I^{(i)}_k$ where $L(x)$ is constant, say $L$. This is described by a sequence of permutations of $\{1,2,\ldots,L\}$, i.e., $\pi_1,\pi_2,\pi_3,\ldots\in S_L$.
For example, for $x\in (\frac{1}{2},\frac{1897}{3762})$, then the combinatorial structure of the M\&m sequence associated to $x$ is as follows:
\begin{align}
x_1&<x_2<x_4<x_5<x_6<x_7<x_9<x_8<x_{10}<x_{17}<x_{18}<x_{13}<x_{14}<x_{21}<x_{15}<x_{22}<\nonumber\\&<x_{16}<x_{29}<x_{30}<x_{31}<x_{32}<x_{37}<x_{38}<x_{39}<x_{40}<x_{41}<x_{42}<x_{11}<x_{45}<x_{46}<\nonumber\\
&<x_{47}<x_{48}<x_{49}<x_{50}<x_{51}<x_{27}<x_
{52}<x_{73}<x_{53}<x_{25}<x_{54}<x_{28}<x_{26}<x_{12}<\nonumber\\
&<x_{23}<x_{24}<x_{33}<x_{59}<x_{60}<x_{61}<x_{62}<x_{63}<x_{34}<x_{
64}<x_{65}<x_{66}<x_{67}<x_{68}<
\nonumber\\
&<x_{69}<x_{70}<x_{71}<x_{72}<x_{19}<x_{20}<x_{35}<x_{36}<x_{57}<x_{58}<x_{55}<x_{43}<x_
{56}<x_{44}<x_{3},\nonumber
\end{align}
(recall that we assume $x_1=1$, $x_2=x$ and $x_3=1$) while for $x\in(\frac{1897}{3762},\frac{919}{1822})$, the last two inequalities 
are replaced by $x_{56}<x_{3}<x_{44}$.
The two permutations $\pi_1=(1,2,4,5,\ldots,43,56,44,3)$ and $\pi_2=(1,2,4,5,\ldots, 43,56,3,44)$ can be obtained from one another by a simple transposition. In general, it is convenient to introduce the permutations $\sigma_1,\sigma_2,\ldots$ satisfying the conditions
$\pi_2=\sigma_1 \pi_1$, $\pi_3=\sigma_2\pi_2$, $\pi_4=\sigma_3\pi_3$, etc. 
For example, for the interval $I^{(0)}_2=(\frac{1}{2},\frac{341}{666}]$ described above we have $L=73$, and we obtain
\begin{align}
\underline{\sigma}^{(I^{(0)}_2)}=&(\sigma_1,\ldots, \sigma_{34}) = ( (72,73), (71,72), (70,71), (69,70), (68,69), (67,68), (66,67), (65,66),\nonumber\\& (64,65),(63,64), (62,63), (61,62), (60,61), (59,60), (58,59), (57,58), (56,57), (55,56),\nonumber\\
&(54,55), (53,54), (52,53), (51,52), (50,51), (49,50), (48,49), (47,48), (46,47), (45,46),\nonumber\\
& (44,45), (43,44), (42,43), (41,42), (40,41), (39,40) ).\nonumber
\end{align}
This can be visualised as in Figure \ref{figure-permutations1} (top).
The transitions $\sigma_i$ are not always so simple. For example, for $I^{(0)}_3=(\frac{341}{666},\frac{24073}{47010}]$ (where $L(x)=75$),
we have
\begin{align}
\underline{\sigma}^{(I^{(0)}_3)}=&(\sigma_1,\ldots, \sigma_{78}) = ((41,42), (40,41), (42,43), (41,42), (43,44), (42,43),
(63,64), (61,62),\nonumber\\
& (62,63), (60,61), (44,45), (43,44),
(60,61), (58,59), (59,60), (57,58), (45,46), (44,45),\nonumber\\
&(58,59), (55,56), (46,47), (45,46), (56,57), (54,55),
(55,56), (53,54), (47,48), (46,47),\nonumber\\
& (54,55), (52,53),
(53,54), (51,52), (50,51), (52,53), (51,52), (49,50),
(50,51), (48,49), \nonumber\\ & (47,48), (49,50), (48,49), (50,51),
(49,50), (51,52), (52,53), (50,51), (51,52), (53,54), \nonumber\\ &
(52,53), (54,55), (46,47), (47,48), (53,54), (55,56),
(54,55), (56,57), (45,46), (46,47), \nonumber\\ &(55,56), (58,59),
(44,45), (45,46), (57,58), (59,60), (58,59), (60,61),
(43,44), (44,45), \nonumber\\& (60,61), (62,63), (61,62), (63,64),
(42,43), (43,44), (41,42), (42,43), (40,41), (41,42)), \nonumber
\end{align}
which is illustrated in Figure \ref{figure-permutations1} (bottom).
Moreover, the $\sigma_i$'s are not always 2-cycles. For example, if $x\in I^{(0)}_{182}=(\frac{706817}{1380274},\frac{1413973}{2761210}]$ (where $L=253$ and $m(x)=\frac{333}{8}x-\frac{325}{16}$) we have
\begin{align}
\underline{\sigma}^{(I^{(0)}_{182})}=&(\sigma_1,\ldots, \sigma_{19})=( (130,131), (129,130), (131,132), (130,131), (132,133),
(131,132), \nonumber\\&(133,134), (132,133), (134,135), (133,134),
(135,136), (134,135), (136,137),\nonumber\\& (135,136), (137,138,139),
(136,137,138), (139,140), (138,139), (140,141)).\nonumber
\end{align}
Product of disjoint cycles can also appear as $\sigma_i$. For example, the first few permutations corresponding to $x$ in the interior of $I^{(1)}_1=I^{(1)}=[\frac{2911001}{5684610},\frac{339}{662}]$ (where $L(x)=271$ and $m(x)=-\frac{2840973}{32}x+\frac{2909701}{64}$) are 
$$\underline{\sigma}^{(I^{(1)}_1)}=(\sigma_1,\ldots,\sigma_{188})=((138,140,139),(139,141)(140,142), (141,143,142),(226,227),\ldots).$$
This analysis shows how the complicated dynamics of the mean-median map is reflected by the combinatorics of the permutations describing M\&m sequences.

\section*{Appendix A}\label{AppendixA}
We collect here all the formul\ae\, found by our algorithm in neighbourhoods of rational points $\frac{1}{2}<\frac{p}{q}<\frac{2}{3}$ with $5\leq q\leq 18$.
Near $\frac{3}{5}$
\begin{align}
m(x)&=
\begin{cases}
-\frac{126466131}{256}x+\frac{75880293}{256},& \frac{76675523}{127792541}\leq x\leq \frac{3}{5}\\
\frac{126310099}{256}x-\frac{75785445}{256},& \frac{3}{5}\leq x\leq \frac{4194392}{6990653}.
\end{cases}\nonumber
\end{align}
Near $\frac{4}{7}$
\begin{align}
m(x)&=
\begin{cases}
-\frac{200985}{16}x+\frac{14359}{2},& \frac{1954924}{3421219}\leq x\leq \frac{4}{7}\\
\frac{202033}{16}x-7214,& \frac{4}{7}\leq x\leq \frac{44485228}{77846827}.\nonumber
\end{cases}\nonumber
\end{align}
Near $\frac{5}{8}$
\begin{align}
m(x)&=
\begin{cases}
-\frac{71307}{4}x+\frac{356603}{32},& \frac{80040341}{128065272}\leq x\leq \frac{5}{8}\\
\frac{71307}{4}x-\frac{356467}{32},& \frac{5}{8}\leq x\leq \frac{149898117}{239835656}.\nonumber
\end{cases}\nonumber
\end{align}
Near $\frac{5}{9}$
\begin{align}
m(x)&=
\begin{cases}
\frac{3}{16} (24185-43521 x),& \frac{24959983}{44928379}\leq x\leq \frac{5}{9}\\
\frac{9}{16} (14467 x-8035),& \frac{5}{9}\leq x\leq \frac{20203}{36365}.
\end{cases}\nonumber
\end{align}
Near $\frac{6}{11}$
\begin{align}
m(x)&=
\begin{cases}
-\frac{139480363}{2048}x+\frac{38041891}{1024},& \frac{2034476}{3729873}\leq x\leq \frac{6}{11}\\
\frac{139446571}{2048}x-\frac{38029091}{1024},& \frac{6}{11}\leq x\leq \frac{2002180}{3670663}.
\end{cases}\nonumber
\end{align}
Near $\frac{7}{11}$
\begin{align}
m(x)&=
\begin{cases}
 -\frac{49767}{16} x+\frac{31691}{16},& \frac{161123}{253207}\leq x\leq \frac{7}{11}\\
\frac{49695}{16} x-\frac{31603}{16} +,& \frac{7}{11}\leq x\leq \frac{1100921}{1729909}.
\end{cases}\nonumber
\end{align}
Near $\frac{7}{12}$
\begin{align}
m(x)&=
\begin{cases}
-213 x+\frac{501}{4},& \frac{7381379}{12670284}\leq x\leq \frac{7}{12}\\
219 x-\frac{507}{4} ,& \frac{7}{12}\leq x\leq \frac{52083521}{89164436}.
\end{cases}\nonumber
\end{align}
Near $\frac{7}{13}$
\begin{align}
m(x)&=
\begin{cases}
- \frac{424918635}{32} x +\frac{228802391}{32} ,& \frac{339327836}{630180267}\leq x\leq \frac{7}{13}\\
\frac{424847139}{32} x-\frac{228763795}{32} ,& \frac{7}{13}\leq x\leq \frac{4675837673}{8683698531}.
\end{cases}\nonumber
\end{align}
Near $\frac{8}{13}$
\begin{align}
m(x)&=
\begin{cases}
- \frac{9504682471614359}{1024} x+\frac{731129420893811}{128}  ,& \frac{411836560641067}{669234411041734}\leq x\leq \frac{8}{13}\\
\frac{9504493870368247}{1024} x-\frac{731114913104851}{128}  ,& \frac{8}{13}\leq x\leq \frac{1656012165640316}{2691019769165513}.
\end{cases}\nonumber
\end{align}
Near $\frac{9}{14}$
\begin{align}
m(x)&=
\begin{cases}
-\frac{313309}{32} x+\frac{402931}{64}  ,& \frac{345573}{537562}\leq x\leq \frac{9}{14}\\
 \frac{312701}{32} x-\frac{401939}{64},& \frac{9}{14}\leq x\leq \frac{83961}{130604}.
\end{cases}\nonumber
\end{align}
Near $\frac{8}{15}$
\begin{align}
m(x)&=
\begin{cases}
-\frac{19695}{8} x+1314,& \frac{15838}{29697}\leq x\leq \frac{8}{15}\\
\frac{19815}{8} x-1320 ,& \frac{8}{15}\leq x\leq \frac{27100}{50811}
\end{cases}\nonumber
\end{align}
Near $\frac{9}{16}$
\begin{align}
m(x)&=
\begin{cases}
- 486 x+\frac{2195}{8},& \frac{599428}{1066971}\leq x\leq \frac{9}{16}\\
494x-\frac{2215}{8},& \frac{9}{16}\leq x\leq \frac{340224257}{604065151}.
\end{cases}\nonumber
\end{align}
Near $\frac{9}{17}$
\begin{align}
m(x)&=
\begin{cases}
-\frac{151759}{32}x+\frac{80375}{32},& \frac{50463303}{95323759}\leq x\leq \frac{9}{17}\\
\frac{152303}{32}x-\frac{80599}{32},& \frac{9}{17}\leq x\leq \frac{82731}{156263}.
\end{cases}\nonumber
\end{align}
Near $\frac{10}{17}$
\begin{align}
m(x)&=
\begin{cases}
-\frac{11203293}{2} x+3295088,& \frac{40009991}{68016985}\leq x\leq \frac{10}{17}\\
\frac{11211513}{2}x-3297502,& \frac{10}{17}\leq x\leq \frac{8901073}{15131824}.
\end{cases}\nonumber
\end{align}
Near $\frac{11}{17}$
\begin{align}
m(x)&=
\begin{cases}
-\frac{11184345087}{1024}x+\frac{7236931237}{1024},& \frac{193668583}{299305992}\leq x\leq \frac{11}{17}\\
\frac{11198078463}{1024}x-\frac{7245813413}{1024},& \frac{11}{17}\leq x\leq \frac{190590017}{294548208}.
\end{cases}\nonumber
\end{align}
Near $\frac{11}{18}$
\begin{align}
m(x)&=
\begin{cases}
-\frac{361818645}{128}x+\frac{442223279}{256},& \frac{362832747}{593726330}\leq x\leq \frac{11}{18}\\
\frac{361521237}{128}x-\frac{441858799}{256},& \frac{11}{18}\leq x\leq \frac{71738975}{117391047}.
\end{cases}\nonumber
\end{align}

\end{document}